     \documentstyle[12pt,fleqn]{article}
      \def\dg{\dagger}
\def\xb{{\bf x}}     \def\yb{{\bf y}}  \def\zb{{\bf z}}

\def\T{\Theta }   \def\a{\alpha }     \def\g{\gamma }
\def\d{\delta }  \def\O{\Omega }  \def\D{\Delta }  \def\e{\epsilon }

\def\beq{\begin{equation}}
\def\eeq{\end{equation}}

\def\iA{A_i(\xb_{\g}) }            \def\jA{A_j(\yb_{\g}) }
\def\kA{A_k(\xb_{\g}) }            \def\iAd{A_i(\xb_{\d}) }
\def\jAd{A_j(\yb_{\d}) }           \def\kAd{A_k(\xb_{\d}) }
     \def\hjA{A_j^{\dg}(\yb_{\g}) }
     
\def\hjAd{A_j^{\dg}(\yb_{\d}) }    
\def\hkAzd{A_k^{\dg}(\zb_{\d}) }

\begin{document}
\begin{titlepage}

\flushright{q-alg/9511017 \\
 ITP-95-24 E }
\vskip 25mm
\begin{center}
{\bf MULTIPARAMETER DEFORMATIONS OF THE ALGEBRA $gl_n$\\
 IN TERMS OF ANYONIC OSCILLATORS}
\end{center}
\vskip 20mm

\begin{center}
{\it A.M.~GAVRILIK,\ \ N.Z.~IORGOV}
\end{center}
\vskip 8mm

\begin{center}
{\it Institute for Theoretical Physics,
National Academy of Sciences of Ukraine\\
Metrologichna str. 14-b, 252143 Kiev, Ukraine}
\end{center}

\vskip 35mm

\begin{abstract}
\noindent Generators of multiparameter deformations
$U_{q;s_1,s_2,...,s_{n-1}}(gl_n)$ of the universal enveloping
algebra $U(gl_n)$ are realized bilinearly by means of appropriately
generalized form of anyonic oscillators (AOs). This modification
takes into account the parameters $s_1,...,s_{n-1}$ and yields usual
AOs when all the $s_i$ are set equal to unity.
\end{abstract}

\end{titlepage}

{\bf 1. Introduction.}
Various aspects of quantum groups and quantum (or $q$-deformed)
algebras \cite{DFRT:,J:} remain to be a subject of intensive study.
Recently, it was shown by Lerda and Sciuto \cite{LS:} that
the $q$-algebra $U_q(su_2)$ admits a realization in terms of
two modes of so-called anyonic oscillators - certain nonlocal
objects defined on a two dimensional square lattice. Shortly after,
this result was extended in \cite{CM:}-\cite{FMS:} to the case of
higher rank $U_q(sl_n)$ algebras and, moreover, to the $q$-analogs
of all semisimple Lie algebras from the classical series
$A_r$,\ $B_r$,\ $C_r$,\ and $D_r$.

During last few years, a {\it multiparameter} (and two-parameter,
in particular) deformations of $GL_n$ groups and of the algebras
$U(gl_n)$, $U(sl_n)$ were developed \cite{DM:}-\cite{S-Q:}.
Accordingly, it is of interest to explore the possibility
of constructing an analogs of anyonic realization for those
multiparametric-deformed algebras. As a step in this direction,
Matheus-Valle and R-Monteiro \cite{MVM:} have presented
a kind of anyonic construction
for the two parametric-deformed algebra $U_{p,q}(sl_2)$.

The subject of our consideration (just in the context of 'anyonic'
realizations) will be a class of multiparameter deformations
$U_{q;s_1,s_2,...,s_{n-1}}(gl_n)$ of the universal enveloping
algebras $U(gl_n)$ which are generated by the elements \
${\bf 1},\ I_{jj+1},\ I_{j+1j}$, \\ $j=1,2,...,n-1$,
and $I_{ii}, \  i=1,2,...,n$,\  and are defined by the relations
\begin{eqnarray}
&&{[I_{ii}\ ,\ \ I_{jj}]} = 0 ,                 \nonumber\\
&&{[I_{ii}\ ,\ \ I_{jj+1}]} = \d_{ij}I_{ij+1} -
\d_{ij+1}I_{ji} ,                                \nonumber\\
&&{[I_{ii}\ ,\ \ I_{j+1j}]} = \d_{ij+1}I_{ij}
- \d_{ij}I_{j+1i} ,                                       \\
&&{[I_{ii+1}, I_{jj+1}]} = {[ I_{i+1i}, I_{j+1j}]} = 0
 \qquad\qquad \hbox{for} \quad \quad \vert {i-j}\vert \ge 2,\nonumber\\
&&{[I_{ii+1}, I_{j+1j}]} = \d_{ij}      \frac    {
   (s_iq)^{I_{ii}}(s_i^{-1}q)^{-I_{i+1,i+1}}
- (s_iq)^{I_{i+1,i+1}}(s_i^{-1}q)^{-I_{ii}}   }
    { s_iq - (s_i^{-1}q)^{-1}  }  \qquad j\ne i\!-\!1, \nonumber\\
&&{[I_{ii+1}, I_{ii-1}]_{(s_i s_{i-1}^{-1}) } } = 0 , \nonumber
\end{eqnarray}
together with the trilinear (or multiparameter generalization of
$q$-Serre) relations

\beq \label{.2}
{\begin{array} {l}
{[\ {[I_{i,i+1}, I_{i+1,i+2}]}_{(s_{i+1}q)},\
    I_{i,i+1}]}_{(s_{i+1}^{-1}q)} = 0,

\vspace{1mm}\\
{[\ {[I_{i+1,i+2}, I_{i,i+1}]}_{(s_{i+1}^{-1}q)},\
    I_{i+1,i+2}]}_{(s_{i+1}q)} = 0,

\vspace{1mm}\\
{[\ {[I_{i+1,i}, I_{i+2,i+1}]}_{(s_i^{-1}q)},\
    I_{i+1,i}]}_{(s_i q)} = 0,

\vspace{1mm}\\
{[\ {[I_{i+2,i+1}, I_{i+1,i}]}_{(s_i q)}, \
    I_{i+2,i+1}]}_{(s_i^{-1}q)} = 0 .
\end{array}       }
\eeq
Here the denotion $[X,Y]_r\equiv XY - r\ YX$ for deformed commutator
is adopted.
With comultiplication rules
\beq \label{.3}
{\begin{array} {l}
\D (I_{kk})\ \ =\ \ I_{kk}\ \otimes\ {\bf 1}\
               +\ {\bf 1}\otimes\ I_{kk},

\vspace{1mm}\\
\D (I_{k,k+1})=I_{k,k+1} \otimes (s_k^{-1}q)^{-I_{kk}}
          \ q^{(I_{kk}+ I_{k+1,k+1})/2}

\vspace{1mm}\\
   \hspace{15mm} \qquad\quad\quad \ + \ (s_k\ q)^{I_{kk}}
   \ q^{-(I_{kk}+ I_{k+1,k+1})/2}\otimes I_{k,k+1},

\vspace{1mm}\\
\D (I_{k+1,k}) = I_{k+1,k} \otimes (s_k\ q)^{I_{k+1,k+1}}
          \ q^{-(I_{kk}+ I_{k+1,k+1})/2}

\vspace{1mm}\\
    \hspace{15mm} \qquad\quad\quad \ + \ (s_k^{-1} q)^{-I_{k+1,k+1}}
  \ q^{(I_{kk}+ I_{k+1,k+1})/2}\otimes I_{k+1,k},
\end{array}        }
\eeq
and the counit structure

\beq  \label{.4}
\e ({\bf 1}) = 1 ,\quad\quad
\e (I_{kk}) =  \e (I_{k,k+1})=
\e (I_{k+1,k}) = 0 ,
\eeq
the considered algebra becomes a bialgebra.  In the
restricted case characterized by

\beq       \label{.5}
s_1 = s_2 = ... = s_{n-1} = s ,
\eeq
the algebra $U_{q;s_1,s_2,...,s_{n-1}}(gl_n)$
reduces to the two-parameter deformation
$U_{q,s}(gl_n)$ given by Takeuchi $\cite{T:}$.
If in addition the restriction $s\!=\!1$ is imposed, the
standard $q$-deformation of Drinfeld--Jimbo $\cite{{DFRT:},{J:}}$ is
recovered. Conversely, the algebra under consideration with defining
relations (1)-(4) can be generated from the standard $U_q(gl_n)$
by applying the procedure described in \cite{R:} (see also
\cite{S-Q:}).

Below, we will present a realization of these
multiparametric-$\!$deformed algebras
$U_{q;s_1,s_2,...,s_{n-1}}(gl_n)$, as defined in
relations (1)-(4), by means of definite set of \  $n$\ modified
anyonic oscillators (or 'quasi-anyons', i.e. an appropriate
generalization of the anyons used in \cite{LS:}-\cite{FMS:}).

\bigskip

{\bf 2. Preliminaries.} Let us begin with some
preliminaries concerning ($d=2$) lattice angle function, disorder
operator, and anyonic oscillators (see \cite{LS:}-\cite{FMS:}).
Let $\O$ be a two dimensional square lattice with the spacing
$a\! =\! 1$.
On this lattice we consider a set of $N$ species (sorts) of fermions
$c_i ({\bf x})$,\ \ $i=1,...,N$,\ \ ${\bf x}\in\Omega$,\ which
satisfy the following standard anti\-commutation relations:
\beq   \label {.6}
\{ c_i ({\bf x}), c_j ({\bf y})\} =
\{ c_i^\dg ({\bf x}), c_j^\dg ({\bf y})\} = 0,
\eeq
\beq         \label {.7}
\{ c_i ({\bf x}), c_j^\dg
({\bf y})\} =\d_{ij}\ \d ({\bf x}, {\bf y}) .
\eeq
Here $\d({\bf x}, {\bf y})$\ is nothing but the conventional
lattice $\d$-function: $\delta ({\bf x}, {\bf y}) = 1$\ if
${\bf x}={\bf y}$\ and vanishes if ${\bf x}\neq {\bf y}$.

We use the same as in \cite{LS:}-\cite{FMS:} definition of the
lattice angle functions $\T_{\g_{\xb }}({\xb },{\yb })$\ and
$\T_{\d_{\xb }}({\xb },{\yb })$\ that correspond to the two
opposite types of cuts ($\g$-type and $\d$-type), and the same
definition of ordering ($\xb > \yb ,\ \xb < \yb $). The
corresponding two types of {\it disorder operators} $K_i ({\bf x}_{\g
})$, and $K_i ({\bf x}_{\d })$, \ $i=1,...,N$, are introduced as
follows:  \beq                \label{.8} {\begin{array} {l} K_j ({\bf
x}_{\g }) = \exp \left ( i\nu \sum_{\yb \neq \xb } \Theta_{\g_\xb
}({\xb }, {\yb })\ c_j^\dg ({\yb }) c_j ({\yb })\right ),

\vspace{1mm}\\

K_j ({\bf x}_{\d }) =
\exp \left ( i\nu \sum_{\yb \neq \xb }
\Theta_{\d_\xb }({\xb }, {\yb })
\ c_j^\dg ({\yb }) c_j ({\yb })\right ).
\end{array}   }
\eeq
The number $\nu $ that appears here is usually called the
{\it statistics parameter}.

The anyonic oscillators (AOs)\ $a_i(\xb_{\g })$ and
$a_i(\xb_{\d })$,\ \ $i=1,...,N$,  are defined \cite{LS:} as
\beq                         \label{.9}
a_i(\xb_{\g }) = K_i ({\bf x}_{\g })\ c_i (\xb ) , \qquad\qquad
a_i(\xb_{\d }) = K_i ({\bf x}_{\d })\ c_i (\xb )
\eeq
(no summation over $i$). One can show that these AOs
satisfy the following relations of permutation. For $i\neq j$ and
arbitrary $\xb ,\yb \in \O $,
\beq                          \label {.10}
\{ a_i (\xb_{\g }), a_j (\yb_{\g })\} =
\{ a_i (\xb_{\g }), a_j^\dg (\yb_{\g })\}=0.
\eeq
Let $q=\exp (i\pi\nu )$. For $i=j$ and for two distinct sites
(i.e. $\xb\ne \yb$) on the lattice $\O$ one has
\beq a_i (\xb_{\g }) a_i         \label {.11}
(\yb_{\g }) +q^{-{\rm sgn }(\xb -\yb )} a_i (\yb_{\g }) a_i (\xb_{\g
}) = 0,
\eeq
\beq a_i (\xb_{\g })a_i^\dg (\yb_{\g }) +q^{{\rm sgn}(\xb -
\yb )}\ a_i^\dg (\yb_{\g }) a_i (\xb_{\g }) = 0, \label{.12}
\eeq
whereas on the same site
\beq \bigl (a_i (\xb_{\g })\ \bigr )^2 = 0,
\qquad\qquad\qquad \{ a_i (\xb_{\g }), a_i^\dg (\xb_{\g })\} = 1.
\label{.13}
\eeq
The analogs of relations (10)-(13) for anyonic oscillators
of the opposite type $\d$ are obtained from (10)-(13) by replacing
$\g\to\d$ and $q\to q^{-1}$.

Note that it is the pair of relations (11), (12) (their analogs for
the $\d$-type of cut, and Hermitian conjugates of all them) which the
statistics parameter $\nu$ does enter. In comparison with ordinary
fermions, the basic feature of anyons is their nonlocality (the
attributed cut) and their braiding property specific of $d=2$ and given
by eq. (11), (12).  These latters imply that anyons of the {\it same
sort}, even allocated at different sites of the lattice, nevertheless
'feel' each other due to the factor expressed by the parameter $q$ (or
$\nu$).

Finally, the commutation relations for anyons of opposite types of
non-locality, i.e. of $\g$-type and $\d$-type, are to be exhibited
($\xb$, $\yb$ arbitrary):
\beq
\{ a_i (\xb_{\g }), a_j (\yb_{\d })\} = 0,    \label{.14}
\eeq
\beq
\{ a_i (\xb_{\g }), a_j^\dg (\yb_{\d })\}=0
\qquad\qquad \hbox{for}\ \ \xb\neq\yb,  \label{.15}
\eeq
\beq \{ a_i (\xb_{\g }), a_j^\dg (\xb_{\d })\} =
\d_{ij} \ q^{\left [\sum_{\yb <\xb } - \sum_{\yb >\xb }\right ] c_i^\dg
(\yb )c_i (\yb )},        \label{.16}
\eeq
as well as the relations that result from all these (i.e. (10)-(16))
by applying Hermitian conjugation.

As proven in \cite{CM:}-\cite{FMS:}, the set of $N$ anyons
defined in (9), through the formulae
\begin{displaymath}
\begin{array} {l}
E_j^+\equiv I_{j,j+1} = \sum_{x\in\O} I_{j,j+1}(\xb ),  \qquad\qquad
E_j^-\equiv I_{j+1,j} = \sum_{x\in\O} I_{j+1,j}(\xb ),

\vspace{2mm} \\
H_j\equiv I_{jj}-I_{j+1,j+1} =
\sum_{x\in\O}\{ I_{jj}(\xb )-I_{j+1,j+1}(\xb )\},
\end{array}
\end{displaymath}
where
\begin{displaymath}
\begin{array} {l}
      I_{j,j+1}(\xb ) = a_j^\dg (\xb_{\g}) a_{j+1}(\xb_{\g}),
\qquad\qquad
I_{j+1,j}(\xb ) = a_{j+1}^\dg (\xb_{\d }) a_j (\xb_{\d }),

\vspace{2mm}\\
I_{jj}(\xb ) = a_j^\dg (\xb_{\a}) a_{j}(\xb_{\a}) = N_j (\xb )
\end{array}
\end{displaymath}
(here $\a$ is $\g$ or $\d$;\ \ $N_j =c_j^\dg c_j$), supplies a
(bilinear) realization of the $U_q(sl_N)$ algebra.

\bigskip

{\bf 3. Quasianyonic oscillators and realization of
$U_{q;s_1,...,s_{n-1}}(gl_n)$.} To realize analogously
the algebra $U_{q;s_1,s_2,...,s_{n-1}}(gl_n)$, we
have to use a modified (with respect to standard definition (9) and to
the relations (11)-(16)) definition of transmuted (from the fermionic
prototypes (6)-(7)) oscillators, namely
\beq \label{.17} {
\begin{array}
{l} \kA = \exp  \left (   i\nu\sum_{\yb\ne \xb}\T_{\g_{\xb}}(\xb, \yb
)N_k(\yb )   \right ) \ \prod_{j=1}^{k-1} s_j^{\sum_{\yb} N_j(\yb )}
c_k(\xb ),

\vspace{2mm} \\

\kAd = \prod_{j=k}^{n-1} s_j^{\sum_{\yb} N_{j+1}(\yb)}
\exp \left ( i\nu\sum_{\yb\ne \xb}\T_{\d_{\xb}}(\xb, \yb) N_k(\yb )
\right )\ c_k(\xb ).
\end{array}   }
\eeq
We'll call them the "quasi-anyonic" operators, since under
restriction (5) and $s\! =\! 1$ they turn into usual anyonic ones,
see eq. (9). Below, let $s_j = exp (i\pi \rho_j ) ,\ \rho_j\in{\bf R} $.

\medskip
By direct examination one verifies the following. For the coinciding
modes of the operators $A_i(\xb_{\a})$ and the same cuts, the relations
of commutation remain the same as those in eqs. (11)-(13). For those
with opposite cuts ($\g$ and $\d$) and at coinciding modes one has
\beq                              \label{.18}
\{ \kA ,\hkAzd \}\!=\d (\xb ,\zb ) q^{ \sum_{\yb\ne \xb} {\rm
sgn} (\xb -\yb ) N_k (\yb ) }\prod_{j=1}^{k-1} s_j^{\sum_{\yb}
N_j(\yb )} \prod_{j'=k}^{n-1} s_{j'}^{- \sum_{\yb} N_{j'+1}(\yb )}.
\eeq
For different modes and the same cut $\g$, we obtain  $\ (i>j)$
\beq                                         \label{.19}
{\begin{array} {l}
\iA \jA + s_j^{-1}\jA \iA = 0 ,

\vspace{2mm}\\

\iA \hjA + s_j\ \hjA \iA = 0 ,
\end{array}   }
\eeq
while for the same cut $\d$ we have $\ (i>j)$
\beq                              \label{.20}
{\begin{array} {l}
\iAd \jAd + s_{i-1} \jAd \iAd = 0 ,

\vspace{2mm}\\
\iAd \hjAd + s_{i-1}^{-1} \hjAd \iAd = 0 .
\end{array}     }
\eeq
Finally, for opposite cuts the relations of permutation are, for
$\forall x, y,$
\beq                                         \label{.21}
\{ \iA , \jAd \}
= 0, \quad i\le j,\quad \{ \iA , \hjAd \} = 0, \quad i < j,\
\eeq
and, again for arbitrary $x, y$ ,
\begin{eqnarray}
\iA \jAd +
s_j^{-1} s_{i-1} \jAd \iA &\!=\!& 0, \quad\quad i>j , \label{.22}\\
\iA \hjAd + s_j s_{i-1}^{-1} \hjAd \iA &\!=\!& 0,
\quad\quad i>j .                  \label {.23}
\end{eqnarray}
To the considered relations for quasianyons, their Hermitian
conjugates are to be added.

It is obvious that at $s_1 = s_2 = .. = s_{n-1} = 1$ all the relations
of permutation for the quasianyonic operators $A_i(\xb_{\a})$ go over
into those for usual anyons.

It can be shown by direct verification that
the following assertion is true.
\medskip

 {\bf Proposition.}\ The generators
$ I_{jj+1}, I_{j+1j}$, \ $j=1,2,...,n-1$,
and $I_{ii}, \ i=1,2,...,n$, realized in the form
\begin{displaymath}
 I_{k,k+1} = \sum_{x\in\O} I_{k,k+1}(\xb ), \qquad
 I_{k+1,k} = \sum_{x\in\O} I_{k+1,k}(\xb ), \qquad
  I_{kk}  = \sum_{x\in\O} I_{kk}(\xb ),
\end{displaymath}
with local densities taken as
\begin{displaymath}
\begin{array} {l}
I_{k,k+1}(\xb ) =
      A_k^\dg (\xb_{\g}) A_{k+1}(\xb_{\g}),     \qquad\qquad
I_{k+1,k}(\xb ) =
      A_{k+1}^\dg (\xb_{\d }) A_k (\xb_{\d }),

\vspace{2mm} \\
I_{kk}(\xb ) =
   A_k^\dg (\xb_{\a}) A_{k}(\xb_{\a}) = N_k (\xb ),
\end{array}
\end{displaymath}
close into the (global bi-)algebra
$U_{q;s_1,s_2,...,s_{n-1}}(gl_n)$ defined by the
relations (1)-(4).

\bigskip
{\bf 4. Concluding remarks.} Let us make some conclusions.
The formulas presented in previous section describe the set
of generalized 'anyons', or quasi-anyons, which possess maximum of
possible inter-mode dependences. Not only the coinciding modes
of our quasi-anyons feel each other (with braiding characterized by
the $q$) at distinct sites of the lattice (this property reproduces
that of usual anyons, see eqns. (11)-(12) and their conjugates)
but moreover, as exhibit relations (19)-(23), the quasi-anyons
participate in {\it graded brading relations} that depend on
the values of indices, thus realizing a kind of ordering
within the set of quasi-anyonic species.

The system of quasi-anyons given by eq. (17) differs from the
modified anyons used in \cite{MVM:}, as shows direct comparison
at $n=2$. Multimode anyon-like deformed oscillators
proposed in refs. \cite{FZ:}, \cite{MP:} (although being not
non-local objects and not tied to specific dimension) resemble
our quasi-anyons of Section 3, and it is useful to analyze
the (dis)similarities in more detail. Finally, there is an
interesting issue concerning the alternative: to attribute
the ($s_j$-dependent) modifying factors in eq. (17) either to
a change of disorder operator (e.g. in the spirit of
ref.\cite{CGFM:}), or to a (multiparameter) deformation of
the multimode fermionic oscillator.

\medskip
This work was partially supported by the International Science
Foundation under the Grant U4J200 and by the Ukrainian State
Foundation for Fundamental Research.

\end{document}